\documentclass[12pt]{amsart}
\newcommand{\proclaim}[2]{\medbreak {\bf #1}{\sl #2} \medbreak}

\def\cal{\mathcal}

\let\newpf\proof \let\proof\relax 
\newenvironment{pf}{\newpf[\proofname]}{\qed\endtrivlist}

\def\0{{\mathbf{0}}}

\def\r{\rho}
\def\tr{\tilde \rho}
\def\g{\gamma}
\def\tg{{\tilde\gamma}}
\def\a{{\tilde a}}
\def\b{{\tilde b}}

\def\d{{\underline d}}

\newcommand{\ntop}[2]{\genfrac{}{}{0pt}{1}{#1}{#2}}

\newtheorem{thm}{Theorem}[section]
\newtheorem{cor}[thm]{Corollary}

\newtheorem{lemma}[thm]{Lemma}

\newtheorem{prop}[thm]{Proposition}

\theoremstyle{remark}
\newtheorem{rem}{Remark}[section]

\numberwithin{equation}{section}

\def \bn {\hfill \\ \smallskip\noindent}

\theoremstyle{definition}

\def\proof{\bn {\bf Proof.} }

\newcommand{\dist}{\operatorname{dist}}

\newcommand{\inter}{\operatorname{int}}

\newcommand{\id}{\operatorname{id}}

\newcommand{\FF}{{\cal F}}

\newcommand{\C}{{\mathbb C}}

\newcommand{\N}{{\mathbb N}}

\newcommand{\R}{{\mathbb R}}

\newcommand{\Z}{{\mathbb Z}}

\def\B0{{\bold{0}}}

\catcode`\@=12

\def\Empty{}
\newcommand\oplabel[1]{
  \def\OpArg{#1} \ifx \OpArg\Empty {} \else
        \label{#1}
  \fi}
                
\newcommand{\comm}[1]{}
\newcommand{\comment}[1]{}

\begin{document}

\title[Quasisymmetric robustness of Collet-Eckmann]
{Quasisymmetric robustness of the Collet-Eckmann condition in the quadratic
family}

\author{Artur Avila and Carlos Gustavo Moreira}

\address{
IMPA -- Estr. D. Castorina 110 \\
22460-320 Rio de Janeiro -- Brazil.
}
\email{avila@impa.br, gugu@impa.br}

\thanks{Partially supported by Faperj and CNPq, Brazil.}

\date{\today}

\begin{abstract}

We consider quasisymmetric reparametrizations of the parameter space
of the quadratic family.
We prove that the set of quadratic maps which are
either regular or Collet-Eckmann with polynomial recurrence
of the critical orbit has full Lebesgue measure.

\end{abstract}

\setcounter{tocdepth}{1}

\maketitle

\tableofcontents

\section{Introduction}

{\it {The intent of this work is just to be a
rigorous reference for \cite {AM1}.  It is based on the
LaTeX file of \cite {AM}: the proofs are very
similar and in many places differ just by change of constants.  In fact the
proofs in \cite {AM} give easily robustness under quasisymmetric
reparametrizations with small qs-constant, there are
some minor differences to remove this assumption (like
introduction of fast landings and bad returns).}}

Here we consider the quadratic family, $f_a=a-x^2$, where $-1/4 \leq a \leq
2$ is the parameter, and we consider its dynamics in the invariant interval.

In \cite {AM}, a thorough understanding of the dynamics of typical quadratic
maps was obtained.  More specifically, it was shown that a typical
quadratic map is either regular (with a periodic attractor)
or Collet-Eckmann (positive Lyapunov exponent of the critical value)
with polynomial recurrence of the critical orbit.  The
first possibility corresponds to a hyperbolic deterministic setting, with
the well known good properties of hyperbolic systems.  The second is a
particularly well studied case of non-uniformly hyperbolic chaotic dynamics:
in the 90's such maps were shown to possess many hyperbolic-like properties
like stochastic stability, exponential decay of correlations and others
(\cite {KN}, \cite {Y2}, \cite {BV} and \cite {BaBeM}).  In particular it
was possible to answer affirmatively Palis Conjecture \cite {Pa}
for the quadratic family.

It was shown in \cite {ALM} that the parameter space of general analytic
families of unimodal maps (with negative Schwarzian derivative) can be
related to the parameter space of quadratic maps through a quasisymmetric
`holonomy map'.  It becomes then feasible to transfer results from the
quadratic family to other families, but there is one obstruction:
quasisymmetric maps are not absolutely continuous.

Here we show that the set of `good' parameters has not only full Lebesgue
measure, but is resistent to a quasisymmetric reparametrization:

\proclaim{Theorem A.}
{
Consider a quasisymmetric reparametrization of the parameter space of the
quadratic family.
The set of parameters which are either regular or Collet-Eckmann:
$$
\liminf_{n \to \infty} \frac {\ln(|Df^n(f(0))|)} {n}>0.
$$
has full Lebesgue measure.
}

\proclaim{Theorem B.}
{
Consider a quasisymmetric reparametrization of the parameter space of the
quadratic family.
The set of parameters which are either regular or
have polynomial recurrence of the critical orbit
$$
0<\limsup_{n \to \infty} \frac {-\ln(f^n(0))} {\ln(n)}<\infty.
$$
has full Lebesgue measure.
} 

In \cite {AM1} this result is used to obtain a proof of Palis Conjecture
for unimodal maps with negative Schwarzian derivative.  Another approach to
this result (with stronger estimates) allows to obtain more general results,
the proof is however more elaborate \cite {A}.

\section{General definitions}

\subsection{Maps of the interval}

Let $f:I \to I$ be a $C^1$ map defined on some interval $I \subset \R$.
The orbit of a point $p \in I$ is the sequence $\{f^k(p)\}_{k=0}^\infty$.
We say that $p$ is recurrent if there exists a subsequence $n_k \to \infty$
such that $\lim f^{n_k}(p)=p$.

We say that $p$ is a periodic point of period $n$ of $f$ if $f^n(p)=p$, and
$n \geq 1$ is minimal with this property.  In this case we say that $p$ is
hyperbolic if $|Df^n(p)|$ is not $0$ or $1$.  Hyperbolic periodic orbits are
attracting or repelling according to $|Df^n(p)|<1$ or $|Df^n(p)|>1$.

We will often consider the distortion of an iterate $f^n$ restricted to
some interval $J \subset I$, such that $f^n|_J$ is a diffeomorphism.  In
this case we will be intrested on the distortion of $f^n|_J$,
$$
\dist(f^n|_J)=\frac {\sup_J Df^n} {\inf_J Df^n}.
$$
This is always a number bigger than or equal to $1$, we will say that it is
small if it is close to $1$.

\subsection{Trees}

We let $\Omega$ denote the set of finite sequences of non zero integers
(including the empty sequence).  Let $\Omega_0$ denote $\Omega$ without the
empty sequence.

We denote $\sigma^+:\Omega_0 \to \Omega$ by
$\sigma^+(j_1,...,j_m)=(j_1,...,j_{m-1})$ and
$\sigma^-:\Omega_0 \to \Omega$ by
$\sigma^-(j_1,...,j_m)=(j_2,...,j_m)$.

\subsection{Growth of functions}

Let $X$ be a class of functions $g:\N \to \R$ such that $\lim_{n \to \infty}
g(n)=\infty$.  We say that a function $f:\N \to \R$ grows at least rate $X$
if there exists a function $g \in X$ such that $f(n) \geq g(n)$ for $n$
sufficiently big.  We say
that it grows at rate $X$ if there are $g_1,g_2 \in X$ such that $g_1(n)
\leq f(n) \leq g_2(n)$ for $n$ sufficiently big.  We say that $f$ decreases
with rate (at least) $X$ if
$1/f$ grows at rate (at least) $X$.

Standard classes are the following.  Linear for linear functions with
positive slope.  Polynomial for functions $g(n)=n^k,k>0$.  Exponential for
functions $g(n)=e^{k n},k>0$.

The standard torrential function $T$ is defined recursively by $T(1)=1$,
$T(n+1)=2^{T(n)}$.  The torrential class is the set of functions
$g(n)=T(\max\{n+k,1\}), k \in \Z$.

Torrential growth can be detected from recurrent estimates easily.
A sufficient condition for a function which is unbounded from above
to grow at least torrentially is an estimate as
$$
f(n+1)>e^{f(n)^a}
$$
for some $a>0$.  Torrential growth is implied by an estimate as
$$
e^{f(n)^a}<f(n+1)<e^{f(n)^b}
$$
with $0<a<b$.

\subsection{Quasisymmetric maps}

Let $\g \geq 1$ be given.
We say that a homeomorphism $f:\R \to \R$ is $\g$-quasisymmetric ($\g$-qs)
if it has a quasiconformal symmetric extension to $\C$
with dilatation bounded by
$\g$.  Notice that quasisymmetric maps form a group under composition and if
$h_1$ is $\g_1$-qs and $h_2$ is $\g_2$-qs then $h_2 \circ h_1$ is $\g_2
\g_1$-qs.

If $X \subset \R$ and $h:X \to \R$ has a $\g$-quasisymmetric extension to
$\R$ we will also say that $h$ is $\g$-qs.

Let $QS(\g)$ be the set of $\g$-qs maps of $\R$.

\section{Real quadratic maps}

If $\lambda \in \C$ we let $f_\lambda:\C \to \C$ denote the (complex)
quadratic map $\lambda-z^2$.  If $\lambda \in \R$ is such
that $-1/4 \leq \lambda \leq 2$ there exists an interval
$I_\lambda=[-\beta,\beta]$ with
$$
\beta=\frac {-1-\sqrt {1+4 \lambda}} {2}
$$
such that $f_\lambda(I_\lambda) \subset I_\lambda$ and $f_\lambda(\partial
I_\lambda) \subset \partial \lambda$.  For such a $\lambda$, the map
$f=f_\lambda|_{I_\lambda}$ is unimodal, that is, is a self map of
$I_\lambda$ with a unique turning point.

We will keep the notation $f_\lambda$ to refer to a quadratic map when we
discuss its complex extension and $f$ to denote a fixed quadratic map when
we discuss its unimodal restriction.

\subsection{The combinatorics of unimodal maps}

In this subsection we fix a real quadratic map $f$ and define some objects
related to it.

\subsubsection{Return maps}

Given an interval $I$ we define the first return map $R_I:X \to I$ where $X
\subset I$ is the set of points $x$ such that there exists $n>0$ with
$f^n(x) \in I$, and $R_I(x)=f^n(x)$ for the minimal $n$ with this property.

\subsubsection{Nice intervals}

An interval $I$ is nice if it is symmetric and the iterates of
$\partial I$ never intersect $\inter I$.  Given a nice interval $I$
we notice that the domain of the
first return map $R_I$ decomposes in a union of intervals
$I^j$, indexed by integer numbers.  We reserve the index $0$ to denote the
component of $0$ if it exists: $0 \in I^0$,
in this case we say $I$ is proper.

If $I$ is nice, it follows that for all $j \in \Z$, $R_I(\partial I^j)
\subset \partial I$.  It follows that if $j \neq 0$ then $R_I|_{I^j}$ is
a diffeomorphism onto $I$ and if $I$ is proper,
$R_I|_{I^0}$ is symmetric with a unique critical point $0$.  As a
consequence, $I^0$ is also a nice interval.

If $R_I(0) \in I^0$, we say that $R_I$ is central.

\subsubsection{Landing maps}

Given a proper interval $I$ we define the landing map $L_I:X \to I^0$
where $X \subset I$ is the set of points $x$ such that there exists
$n \geq 0$ with
$f^n(x) \in I^0$, and $L_I(x)=f^n(x)$ for the minimal
$n$ with this property.  We notice that $L_I|_{I^0}=\id$.

\subsubsection{Trees}

If $I$ is a proper interval, the first return map to $I$ naturally relates
to the first landing to $I^0$ in the following way.

If $\d \in \Omega$, we define
$I^{\d}$ inductively in the following way.  $I^{\d}=I$
if $\d$ is empty and if $\d=(j_1,...,j_m)$ we let
$I^{\d}=(R_I|_{I_{j_1}})^{-1}(I^{\sigma^-(\d)})$.

We denote $R^{\d}_I=R_I^{|\d|}|_{I^{\d}}$
which is always a diffeomorphism over $I$.

We denote by $C^{\d}=(R^{\d}_I)^{-1}(I^0)$.  We notice
that the domain of the first landing map $L_I$ coincides with the
union of the $C^{\d}$, and furthermore $L_I|_{C^{\d}}=R^{\d}_I$.

Notice that this allows us to relate $R_I$ to $R_{I^0}$, since
$R_{I^0}=L_I \circ R_I$.

\subsubsection{Renormalization}

We say that $f$ is renormalizable if there is an interval $I$ with
and $m>1$ such that $f^m(T) \subset T$ and
$f^j(T) \cap I=\emptyset$ for $1 \leq j<n$.
The maximal such interval is called the renormalization interval of
period $m$, it has the property that $f^m(\partial T) \subset \partial T$.

The set of renormalization periods of $f$ gives an increasing (possibly
empty) sequence of numbers $m_i$, each related to a unique renormalization
interval $T^{(i)}$ which form a nested sequence of intervals.

We say that $f$ is finitely renormalizable if there is a smallest
renormalization interval $T^{(k)}$.  We say that
$f \in \FF$ if $f$ is finitely
renormalizable and $0$ is recurrent but not periodic.

\subsubsection{Principal nest}

Let $\Delta_k$ denote the set of all maps $f$
which have (at least) $k$ renormalizations and which have an orientation
reversing non-attracting periodic point of period $m_k$ which we denote
$p_k$ (that is, $f^{m_k}|_{T^{(k)}}$ has a fixed point
$p_k$ and $Df^n(p_k) \leq -1$).  In this case we
denote by $T^{(k)}_1=[-p_k,p_k]$.  We define by induction a (possibly
finite) sequence $T^{(k)}_i$, such that $T^{(k)}_{i+1}$ is the component of
the domain of $R_{T^{(k)}_i}$ containing $0$.  If this sequence is
infinite, then either it converges to a point or to an interval.

In the former case, $f$ has a recurrent critical point which is not
periodic, and it is possible to show that $f$ is not $k+1$ times
renormalizable, obviously in this case we have $f \in \FF$, and we write
$\FF_k$ for the set of such maps.  If the limit set is an
interval, it is possible to show that $f$ is $k+1$ times renormalizable.

We can of course write $\FF$ as a disjoint union
$\cup_{i=0}^\infty \FF_i$.  For a map
$f \in \FF_k$ we refer to the sequence
$\{T^{(k)}_i\}_{i=1}^\infty$ as the principal nest.

It is important to notice that the domain of the first return map
to $T^{(k)}_i$ is always dense in $T^{(k)}_i$.  Moreover, outside
a very special case, the return map has a hyperbolic structure.

\begin{lemma} \label {hyperbol}

Assume $T^{(k)}_i$ does not have a non-hyperbolic periodic orbit in its
boundary.  For all $T^{(k)}_i$ there exists $C>0$, $\lambda>1$
such that if $x,f(x),...,f^{n-1}(x)$ do not belong to
$T^{(k)}_i$ then $|Df^n(x)|>C \lambda^n$.

\end{lemma}

This theorem is a simple consequence of a general theorem of Guckenheimer
on hyperbolicity of maps of the interval without critical points and
non-hyperbolic orbits (Guckenheimer assumes negative Schwarzian derivative,
so this applies directly to our case, the general case is also true by
Ma\~n\'e's Theorem, see \cite {MvS}).  Notice that the existence or a
non-hyperbolic periodic orbit in the boundary of $T^{(k)}_i$ depends on a
very special combinatorial setting, in particular, all $T^{(k)}_j$ must
coincide (with $[-p_k,p_k]$), and the $k$-th renormalization of $f$ is in
fact renormalizable of period $2$.

This Lemma shows that the complement of the domains of the return map to
$T^{(k)}_i$ (outside the special case) form a regular Cantor set.  This has
many useful consequences (for instance, the image of a regular Cantor set by
a quasisymmetric map has always $0$ Lebesgue measure.

\subsubsection{Lyubich's Regular or Stochastic dichotomy}

A map $f \in \FF_k$ is called simple if the principal nest has
only finitely many central returns, that is, there are
only finitely many $j$ such that $R|_{T^{(k)}_j}$ is central.

In \cite {parapuzzle}, it was proved that almost every
quadratic map is either regular or simple or infinitely renormalizable.  It
was then shown in \cite {regular} that infinitely renormalizable maps have
$0$ Lebesgue measure.

In \cite {ALM} it is remarked that Lyubich's estimates actually prove
quasisymmetric robustness of the set of regular or simple maps, that is,
they are still typical after reparametrization.

\subsubsection{Strategy}

Both Theorems A and B will be proved using the same strategy.

Our strategy is to describe the dynamics of the principal nest, this is our
phase analysis.  From time to time, we transfer the information from the
phase space to the parameter, following the description of the parapuzzle
nest which we will make in the next subsection.

Due to Lyubich's results, we can completely forget about infinitely
renormalizable maps, we just have to prove the claimed estimates for almost
every (after reparametrization) simple map.  During our discussion,
for notational reasons, we will
actually consider a fixed renormalization level $k$, that is, we will
analyse maps in $\Delta_k$.

This allow us to fix some convenient notation: given $g \in \Delta_k$
we define $I_i[g]=T^{(k)}_i[g]$, so that $\{I_i[g]\}$
is a sequence of intervals (possibly finite).
We use the notation $R_i[g]=R_{I_i[g]}$, $L_i[g]=L_{I_i[g]}$ and so on. 
When doing phase analysis (working with fixed $f$) we usually drop the
dependence on the map and write $R_i$ for $R_i[f]$.

\subsection{Parameter partition}

Part of our work is to transfer information from the phase space of some
map $f \in \FF$ to a neighborhood of its parameter space.  This is done in
the following way.  We consider the first landing map to $I_i$.
The domains of this map partition the interval $I_i$, the complement of
this set is a hyperbolic Cantor set $K_i=I_i \setminus \cup C^\d_i$.
This Cantor set persists in a small parameter
neighborhood $J_i$ of $f$, changing in a
continuous way.

Along $J_i$, the first landing map is topologically the same (in a way
will be clear soon).  However the critical value $R_i[g](0)$ moves relative
to the partition (when $g$ moves in $J_i$).  This allows us to partition
the parameter piece $J_i$ in smaller pieces, each
corresponding to a region where $R_i(0)$ belongs to some fixed domain of the
first landing map.

For a discussion of the next Theorem, see \cite {AM}.

\proclaim{Topological Phase-Parameter relation.}
{
Let $f \in \FF$ be exactly $k$ times renormalizable.
Then there is a sequence $\{J_i\}_{i \in \N}$ of nested
parameter intervals (the principal parapuzzle nest of $f$)
with the following properties.

\begin{enumerate}

\item $J_i$ is the maximal interval containing $f$ such that for all $g \in
J_i$ the interval $T^{(k)}_{i+1}[g]$ is defined and changes in a
continuous way.  Since the first return map to $R_i[g]$ has a central
domain, so that the landing map
$L_i[g]:\cup C^\d_i[g] \to I_i[g]$ is defined.

\item $L_i[g]$ is topologically the same along $J_i$: there exists
homeomorphisms $H_i[g]:I_i \to I_i[g]$, such that
$H_i[g](C^\d_i)=C^\d_i[g]$.
The maps $H_i[g]$ change continuously.

\item There exists a homeomorphism $\Xi_i:I_i \to J_i$ such that
$\Xi_i(C^\d_i)$ is
the set of $g$ such that $R_i[g](0)$ belongs to $C^\d_i[g]$.

\end{enumerate}
}

The homeomorphisms $H_i$ and $\Xi_i$ are not uniquely defined, it is easy to
see that we can modify then inside each $C^\d_i$ window.  $H_i$ and $\Xi_i$
are well defined maps if restricted to $K_i$.

With this result we can define for any $f \in \FF$ intervals
$J^j_i=\Xi_i(I^j_i)$ and $J^\d_i=\Xi_i(I^\d_i)$.  From the description we
gave it immediately follows that two intervals $J_i[f]$ and $J_i[g]$
associated to different
maps $f$ and $g$ are either disjoint or nested, and the same happens for
intervals $J^j_i$ or $J^\d_i$.

We will concentrate on the analysis of the regularity of $\Xi_i$ for the
special class of simple maps $f$: one of the good properties of the class of
simple maps is better control of the phase-parameter relation.
Even for simple maps, however, the regularity of $\Xi_i$ is not great:
there is too much dynamical information contained in it.  A solution to this
problem is to forget some dynamical information.  With this intent we
introduce an interval which will be used to erase information.

\subsubsection{Gape interval}

If $i>1$, we define the gape interval $\tilde I_{i+1}$ as follows.

We have that
$R_i|_{I_{i+1}}=L_{i-1} \circ R_{i-1}=R^{\d}_{i-1} \circ R_{i-1}$
for some $\d$, so that
$I_{i+1}=(R_{i-1}|_{I_i})^{-1}(C^{\d}_{i-1})$.  We define the
gape interval
$\tilde I_{i+1}=(R_{i-1}|_{I_i})^{-1}(I^{\d}_{i-1})$.

We notice that for each $I^j_i$, the gape interval
$\tilde I_{i+1}$ either contains or is disjoint from $I^j_i$.

\subsubsection{The phase-parameter relation}

As we discussed before, the dynamical information contained in $\Xi_i$
is entirely given by $\Xi_i|_{K_i}$: a map obtained by $\Xi_i$ by
modification inside a $C^\d_i$ window has still the same properties. 
Therefore it makes sense to ask about the regularity of $\Xi_i|_{K_i}$.  As
we anticipated before we must erase some information to obtain good results. 

Let $f$ be a simple map and $\tau_i$ be such that
$R_i(0) \in I^{\tau_i}_i$.  We define two Cantor sets,
$K^\tau_i=K_i \cap I^{\tau_i}_i$ which
contains refined information restricted to the $I^{\tau_i}_i$ window and
$\tilde K_i=I_i \setminus \cup I^j_i \setminus \tilde I_{i+1}$, which
contains global information, at the cost of erasing information inside each
$I^j_i$ window and in $\tilde I_{i+1}$.

The following Theorem is one of the main steps in \cite {AM}.

\proclaim{Phase-Parameter relation.}
{
Let $f$ be a simple map.  For all $\g>1$ there exists $i_0$ such that
for all $i>i_0$ we have

\begin{description}

\item[Quadratic PhPa1] $\Xi_i|_{K^\tau_i}$ is $\g$-qs,

\item[Quadratic PhPa2] $\Xi_i|_{\tilde K_i}$ is $\g$-qs,

\item[PaPa1] $H_i[g]|_{K_i}$ is $\g$-qs if $g \in J^{\tau_i}_i$,

\item[PaPa2] the map $H_i[g]|_{\tilde K_i}$ is
$\g$-qs if $g \in J_i$.

\end{description}
}

\subsection{Reparametrization}

To prove Theorems A and B, we have to show that for any $\g_0$, if $h$
is $\g_0$-qs then almost every quadratic map
{\it {after reparametrization by $h$}} have some good properties. 
Therefore, we fix now such $\g_0$ and such a reparametrization.  We keep
otherwise the same notation:

\proclaim{Reparametrized Phase-Parameter relation.}
{
Let $f$ be a simple map.  For all $\lambda>1$ there exists $i_0$ such that
for all $i>i_0$ we have

\begin{description}

\item[PhPa1] $h \circ \Xi_i|_{K^\tau_i}$ is $\lambda\g_0$-qs,

\item[PhPa2] $h \circ \Xi_i|_{\tilde K_i}$ is $\lambda\g_0$-qs,

\item[PaPa1] $H_i[g]|_{K_i}$ is $\lambda$-qs if $g \in J^{\tau_i}_i$,

\item[PaPa2] the map $H_i[g]|_{\tilde K_i}$ is
$\lambda$-qs if $g \in J_i$.

\end{description}
}

Notice that PaPa1 and PaPa2 did not change after reparametrization, since
they are just Phase-Phase estimates.

{\bf {All references to the parameter space from now on will assume this
reparametrization.}}

\section{Measure and capacities} \label {measure}

\subsection{Quasisymmetric maps}

If $X \subset \R$ is measurable, let's denote $|X|$ its Lebesgue measure.
Let's describe more metric properties of $\g$-qs maps.

To each $\g$, there exists a constant $k \geq 1$ such that for all
$f \in QS(\g)$, for all $J \subset I$ intervals,
$$
\frac {1} {k} \left ( \frac {|J|} {|I|} \right )^k \leq \frac {|f(J)|}
{|f(I)|} \leq
\left ( \frac {k|J|} {|I|} \right )^{1/k}.
$$

Furthermore $\lim_{\g \to 1} k(\g)=1$.
So for each $\epsilon>0$ there exists $\g>1$ such that $k(2
\g-1)<\epsilon/3$.

\subsection{Capacities and trees}

The $\g$-capacity of a subset $X$ of an interval $I$ is defined as
follows.
$$
p_\g(X|I)=\sup_{h \in QS(\g)} \frac {|h(X)|} {|h(I)|}.
$$

This geometric quantity is well adapted to our context, since it is well
behaved under tree decompositions of sets.  In other words, if $I^j$ are
disjoint subintervals of $I$ and $X \subset \cup I_j$ then
$$
p_\g(X|I) \leq p_\g(\cup I_j|I) \sup p_\g(X \cap I^j|I^j).
$$

\subsection{A measure theoretical lemma}

Our usual procedure consists in picking a class $X$ of maps which we show is
full measure among non-regular maps and then for each map in this class
we describe what happens
for the principal nest, showing finally that a subset $Y$ of $X$
is still full measure.  The first such step was to consider the class of
simple maps.

We describe here our usual argument (a variation of Borel-Cantelli lemma).
Assume at some point we know how to prove that almost every map
belong to some set $X$.  Let $Q_n$ be a (bad) property that a map
may have (usually related to the $n$-th stage of the principle nest).
Suppose we prove that if $f \in X$ then
the probability that a map in $J_n(f)$ has the property $Q_n$ is bounded by
$q_n(f)$ which is shown to be summable
for all $f \in X$.  We then conclude that almost every map does not have
property $Q_n$ for $n$ big enough.

Sometimes we also apply the same argument, proving instead that $q_n(f)$
is summable  where $q_n(f)$ is the
probability that a map in $J^{\tau_n}_n(f)$ has property $Q_n$, where
$\tau_n$ is such that $f \in J^{\tau_n}_n(f)$.

In other words, we use the following lemma.

\begin{lemma}

Let $X \subset \R$ be a measurable set such that for each
$x \in X$ is defined a sequence
$J_i(x)$ of nested intervals converging to $x$
such that for all $x_1,x_2 \in X$
and any $i$, $J_i(x_1)$ is either equal or disjoint to $J_i(x_2)$.  Let
$Q_n$ be measurable subsets of
$\R$ and $q_n(x)=|Q_n \cap J_n(x)|/|J_n(x)|$.  Let $Y$ be
the set of $x$ in $X$ which belong to finitely many $Q_n$.
If $\sum q_n(x)$ is finite for almost any $x \in X$ then $|Y|=|X|$.

\end{lemma}

\begin{pf}

Let $Y_n=\{x \in X|\sum_{k=n}^\infty q_k(x)<1/2\}$.  It is clear that $Y_n
\subset Y_{n+1}$ and $|\cup Y_n|=|X|$.

Let $Z_n=\{x \in Y_n||Y_n \cap J_m(x)|/|J_m(x)|>1/2, m \geq n\}$.
It is clear that $Z_n \subset Z_{n+1}$ and $|\cup Z_n|=|X|$.

Let $T^m_n=\cup_{x \in Z_n} J_m(x)$.  Let $K^m_n=T^m_n \cap Q_m$.
Of course
$$
|K^m_n| \leq \int_{T^m_n} q_m \leq 2 \int_{Y_n} q_m.
$$

And of course
$$
\sum_{m \geq n} \int_{Y_n} q_m \leq \frac{1}{2} |Y_n|.
$$

This shows that $\sum_{m \geq n} |K^m_n| \leq |Y_n|$,
so almost every point in $Z_n$ belong to finitely many $K^m_n$.
We conclude then that almost every point in
$X$ belong to finitely many $Q_m$.

\end{pf}

In practice, we will estimate the capacity of sets in the phase space:
that is, given a map $f$ we will obtain subsets $\tilde Q_n(f)$ in the phase
space, corresponding to bad branches of return or landing maps.  We will the
show that for some $\g>\g_0$ we have
$\sum p_\g(\tilde Q_n(f)|J_n(f))<\infty$ or
$\sum p_\g(\tilde Q_n(f)|J^{\tau_n}_n(f))<\infty$.  We will then use PhPa1,
PhPa2 and the measure-theoretical lemma to conclude that with total
probability among non-regular maps, for all $n$ sufficiently big,
$R_n(0)$ does not belong to a bad set.

From now on when we prove that almost every non-regular parameter
has some property, we will just say that with total probability
(without specifying) such property holds.

\section{Statistics of the principal nest}

Through the end of this paper we fix some constant $\g>\g_0$.  We also fix
$\b \gg 1000 k(2 \g-1)^1000$ and $\a=b^{-1}$, and set $b=\b^{1000\b}$,
$a=b^{-1}$.

\subsection{Decay of geometry}

Let as before $\tau_n \in \Z$ such that $R_n(0) \in I^{\tau_n}_n$.

An important parameter in our construction will be the scaling factor
$$
c_n=\frac {|I_{n+1}|} {|I_n|}.
$$
This variable of course changes inside each
$J^{\tau_n}_n$ window, however, not by much.
From PaPa1, for instance, we get
that with total probability
$$
\lim_{n \to \infty} \sup_{g_1,g_2 \in J^{\tau_n}_n} \frac {\ln(c_n[g_1])}
{\ln(c_n[g_2])}=1.
$$

This variable is by far the most important on our analysis of the statistics
of return maps.  We will often consider other variables
(say, return times): we will show that the distribution of those
variables is concentrated near some average value.
Our estimates will usually give a range
of values near the average, and $c_n$ will play an important role.  Due
(among other issues) to the variability of $c_n$ inside the parameter
windows, the ranges we select will depend on $c_n$ up to an exponent
(say, in a range from $a$ to $b$).  From the estimate we just obtained,
for big $n$ the variability (margin of error) of $c_n$ will fall
confortably in such range, and we won't elaborate more.

If $x \in I^j_n$ we let $j^{(n)}(x)=j$ and if
$x \in C^\d_n$ we let $\d^{(n)}(x)=\d$.

\begin{lemma} \label {estimate on m}
  
With total probability, for all $n$ sufficiently big we have
$$
p_{2 \g-1}(|\d^{(n)}(x)| \leq k|x \in I_n)<k c_n^\a,
$$
$$
p_{2 \g-1}(|\d^{(n)}(x)| \geq k|x \in I_n)<e^{-k c_n^\b}.
$$
  
We also have
$$
p_{2 \g-1}(|\d^{(n)}(x)| \leq k|x \in I^{\tau_n}_n)<k c_n^\a,
$$
$$
p_{2 \g-1}(|\d^{(n)}(x)| \geq k|x \in I^{\tau_n}_n)<e^{-k c_n^\b}.
$$
  
\end{lemma}

\begin{pf}

Let's compute the first two estimates.

Since $I^0_n$ is in the middle of
$I_n$, we have as a simple consequence of the real Schwarz Lemma (see \cite
{attractors}) that
$$
\frac {c_n}{4}<\frac {|C^\d_n|}{|I^\d_n|}<4 c_n.
$$
As a consequence
$$
p_{2 \g-1}(|\d^{(n)}(x)|=m|x \in I_n)<(4 c_n)^{2 \a}
$$
and we get the estimate summing up on $0 \leq m \leq k$.

For the same reason, we get that
$$
p_{2 \g-1}(|\d^{(n)}(x)>m|x \in I_n)<1-(c_n/4)^{\b/2}
p_{2 \g-1}(|\d^{(n)}(x)| \geq m|x \in I_n).
$$
This implies
$$
p_{2 \g-1}(|\d(x)| \geq m|x \in I_n) \leq (1-(c_n/4)^{\b/2})^m .
$$
  
The second estimate follows from
$$
(1-(\frac {c_n}{4})^{\b/2})^k<
(1-c_n^\b)^k<((1-c_n^b)^{c_n^{-\b}})^{k
c_n^\b}<e^{-k c_n^\b}.
$$ 
   
The two remaining estimates are analogous.
\end{pf}

Transferring the result (more precisely the second pair of estimates)
to the parameter in each
$J^{\tau_n}_n$ window using PhPa1 we get (noticing
that the measure of the
complement of the set of parameters in $J^{\tau_n}_n$ such that
$c_n^{-\a/2}<s_n<c_n^{-2\b}$ can be estimated by
$2 c_n^{\a/2}$ for $n$ big, which is summable.

\begin{lemma} \label {growth of s_n}

With total probability,
$$
c_n^{-\a/2}<\frac {\ln (s_n)} {\ln (c^{-1}_n)}<c_n^{-2\b}.
$$
\end{lemma}

\begin{rem}

The parameter $s_n$ influences the size of $c_{n+1}$ in a determinant way. 
It is easy to see (using for instance the real Schwarz Lemma, see \cite
{attractors}) for instance that $\ln(c_{n+1}^{-1})>K s_n$ for some
universal constant $K$, which in general is bounded from below (real a
priori bounds), but since we have decay of geometry, actually goes to
infinity.

\end{rem}

As an easy consequence we get

\begin{cor} \label {c_n torrential}

With total probability,
$$
\frac {\ln(\ln(c_{n+1}^{-1}))} {\ln(c_n^{-1})} > \a/3.
$$
In particular, $c_n$
decreases at least torrentially fast.

\end{cor}

\subsection{Fine partitions}

We use Cantor set $K_n$ and $\tilde K_n$ to partition the phase space.  In
many circumstances we are directly concerned with intervals of this
partition.  However, sometimes we just want to exclude an interval of given
size (usually around $0$).  This size does not usually correspond to a union
of gaps, so we instead should substitute in applications an interval which
is union of gaps, with approximately the given size.  The degree of relative
approximation will always be torrentially good (in $n$), so we usually won't
elaborate on this.  In this section we just give some results which will
imply that the partition induced by the Cantor sets are fine enough to allow
torrentially good approximations.

The following Lemma summarizes the situation.  The proof is based on
estimates of distortion using the real Schwarz Lemma and
the Koebe Principle (see \cite {attractors}) and is very simple,
so we skip it.

\begin{lemma}

The following estimates hold.

\begin{itemize}

\item $|I^j_n|/|I_n|=O(\sqrt {c_{n-2}})$,

\item $|I^{\d}_n|/|I^{\sigma^+(\d)}_n|=O(\sqrt
{c_{n-2}})$.

\item $|C^{\d}_n|/|I^{\d}_n| \approx c_n$,

\item $|\tilde I_{n+1}|/|I_n|=O(e^{1-s_{n-1}})$.

\end{itemize}

\end{lemma}

In other words, distances in $I_n$ can be measured with precision $\sqrt
{c_{n-2}}|I_n|$ in the partition induced by $\tilde K_n$, due to first and
last items (since $e^{1-s_{n-1}}=O(c_{n-1})$).

Distances can be measured much more precisely with respect to the partition
induced by $K_n$, in fact we have good precision in each
$I^{\d}_n$ scale.  In other words, inside $I^{\d}_n$,
the central gap $C^{\d}_n$ is
of size $O(c_n|I^{\d}_n|)$ (third item)
and the other gaps have size $O(\sqrt {c_{n-2}} |C^{\d}_n|)$
(second and third item).

\begin{rem}

We need to consider intervals which are union of gaps due to our phrasing of
the phase-parameter relation, which only gives information about such gaps. 
However, this is not absolutely necessary, and we could have proceeded in a
different way: our proof of the phase-parameter relation actually shows that
there is a holonomy map with good qs estimates between phase and parameter
intervals (and not only Cantor sets).  While this map is not canonical,
the fact that it is a holonomy map for a motion with good phase-phase
estimates would allow our proofs to work.

\end{rem}

\subsection{Initial estimates on distortion}

To deal with the distortion control we need some preliminary known results.
We won't get too much in details here, those estimates are related to the
estimates on gaps of Cantor sets and the Koebe Principle, and can be
concluded easily.

\begin{prop}

For any $j$, if $R_n|_{I^j_n}=f^k$ then
$\dist(f^{k-1}|_{f(I^j_n)})=1+O(c_{n-1})$.

For any $\d$,
$\dist(R^{\sigma^+(\d)}_n|_{I^{\d}_n})=1+O(\sqrt
{c_{n-2}})$.

\end{prop}

We will use the following immediate consequence for the decomposition of
certain branches.

\begin{lemma} \label{decomposition}

With total probability,

\begin{itemize}

\item $R_n|_{I^0_n}=\phi \circ f$ where $\phi$ has torrentially small
distortion,

\item $R^{\d}_n=\phi_2 \circ f \circ \phi_1$ where $\phi_2$ and
$\phi_1$ have torrentially small distortion and
$\phi_1=R^{\sigma^+(\d)}_n$.

\end{itemize}

\end{lemma}
\subsection{Estimating derivatives}

\begin{lemma} \label {away from the boundary}

With total probability,
$$
d(R_n(0),\partial I_n \cup \{0\})<|I_n|n^{-\b}.
$$
In particular $R_n(0) \notin \tilde I_{n+1}$ for all $n$ large enough.

\end{lemma}

\begin{pf}

This is a simple consequence of PhPa2, using that
$n^{-\sqrt \b}$ is summable. 
\end{pf}

From now on we suppose that $f$ satisfies the conditions of the above lemma.

\begin{lemma} \label {dist}

With total probability, for $n$ big enough and $j \neq 0$
$$
\dist(f|_{I^j_n}) \leq n^\b.
$$

\end{lemma}

\begin{pf}

Notice that the gaps of the Cantor sets $K_n$
inside $I^\d_n$ which are different from $C^\d_n$
are torrentially (in $n$) smaller then $C^\d_n$.

Denote by $P^{\d}_n$ a $|C^{\d}_n|/n^\b$
neighborhood of $C^{\d}_n$.

It is clear that the image by a $\g$-qs homeomorphism
of $P^\d_n \setminus C^\d_n$
is of order $n^{-\sqrt \b}|C^\d_n|$
(this is true both on $I_n$ or restricting to
$I^{\tau_n}_n$).  Since $C^\d_n$ are disjoint,
$$
p_\g(I^{\tau_n}_n \cap \cup (P^\d_n \setminus C^\d_n)|I^{\tau_n}_n)
$$
is summable.

Transferring the result to the parameter using
PhPa1 we see that the critical point will
never be in a $n^{-\b}$ neighborhood of any $I^j_{n+1}$ with
$j \neq 0$.
\end{pf}

Applying Lemma \ref {decomposition} we get

\begin{lemma} \label {distortion}

With total probability, for $n$ big enough and for all $\d$
$$
\dist(R^{\d}_n) \leq n^\b.
$$
In particular, for $n$ big enough,
$\sup_{\d \neq \emptyset} \dist(R^{\d}_n)$
is bounded by $2^n$ and
$\inf |R'_n|_{\cup_{j \neq 0} I^j_n}|>2$.

\end{lemma}

\begin{rem} \label {approx}

Lemma \ref {dist} has also an application for approximation of
intervals.  The result implies that if $I^j_n=(c,d)$ and $j \neq 0$, we have
$1/2^n<d/c<2^n$.  As a consequence, for all symmetric interval $I_{n+1}
\subset X \subset I_n$, there exists a symmetric interval $X \subset \tilde
X$, which is union of $I^j_n$ and such that $|\tilde X|/|X|<2^n$
(approximation by union of $C^\d_n$, with $|tilde X|/|X|$ torrentially close
to $1$, follows more easily from the discussion on fine partitions).

\end{rem}

We will also need to estimate derivatives of iterates of $f$, and not only
of return branches.

\begin{lemma} \label {lower bound}

With total probability, if $n$ is sufficiently big and if
$x \in \cup_{j \neq 0} I^j_n$ and $R_n|_{I^j_n}=f^r$, then for $1 \leq k
\leq r$, $|(Df^k(x))|>|x| c_{n-1}^3$.

\end{lemma}

\begin{pf}

Let $n_0$ be minimum such that if $n \geq n_0$ then
$$
\inf |R'_n|_{\cup_{j \neq 0} I^j_n}|>1.
$$
From hyperbolicity of $f$, from Lemma \ref {hyperbol},
restricted to the complement of
$I_{n_0}$, there exists a constant $C>0$ such that if
$f^s(x) \notin I^0_{n_0}$, $r \leq s \leq k$ then $|Df^{k-r}(f^r(x))|>C$.

If $k=r$, the result follows from Lemma \ref {distortion}.  We assume $k<r$.

Let's define $d(s)$, $1 \leq s \leq k$ such that $f^s(x) \in I_{d(s)}
\setminus I^0_{d(s)}$. 
Let $m(s)=\max_{s \leq t \leq k} d(s)$.
Let $k_0=0$.
Supposing $k_j<k$ define $k_{j+1}$ as $k_{j+1}=\max \{k_j<s \leq k|
d(s)=m_s\}$.

If $d(k_1)=n$, $k=r$ and from Lemma \ref {distortion} $|Df^k(x)|>1$.

We have then a sequence $0=k_0 < k_1<...<k_l=k$.
Let $k_i$ be maximal with $d(k_i) \geq n_0$.  We have
$$
|(Df^{k-k_i}(f^{k_i}(x))|>C|f^{k_i}(x)|
$$
and for $1< j \leq i$,
$$
|Df^{k_j-k_{j-1}}(f^{k_{j-1}}(x))|>|Df(f^{k_{j-1}})(x)|>c_{d(k_{j-1})}.
$$
We also have $|Df^{k_1}(x)| \geq |x|$ if $i \geq 1$.  If $i=0$,
$|Df^k(x)| \leq C|x|$.  Combining it all we get
$$
|Df^k(x)|>C |x| \prod_{1 \leq j \leq i} c_{d(k_j)}>|x| c_{n-1}^3.
$$
\end{pf}

\section{Sequence of quasisymmetric constants and trees}

\subsection{Preliminary estimates}

We will need from now on to consider not only $\g$-capacities with $\gamma$
fixed, but a sequence $\g_n$ converging to $\g$.

We define the sequences $\r_n=(n+1)/n$ and $\tr_n=(2n+3)/(2n+1)$,
so that $\r_n>\tr_n>\r_{n+1}$ and $\lim \r_n=1$.
We define the sequence $\g_n=\g \r_n$ and an intermediate sequence
$\tg_n=\g \tr_n$.

As we know, renormalization proccess has two phases, first
$R_n$ to $L_n$ and then $L_n$ to $R_{n+1}$.  The following
remarks shows why it is useful to consider the sequence
of quasisymmetric constants due to losses related to distortion.

\begin{rem} \label {remark 1}

Let $S$ be an interval contained in $I^{\d}_n$.
Using Lemma \ref {decomposition}
we have $R^\d_n|_S=\psi_2 \circ f \circ \psi_1$,
where the distortion of $\psi_2$ and $\psi_1$ are torrentially small
and $\psi_1(S)$ is contained in some $I^j_n$, $j \neq 0$.
If $S$ is contained in $I^0_n$ we may as well write $R_n|_S=\phi \circ f$,
with $\dist(\phi)$ torrentially small.

In either case, if we decompose $S$ in $2km$ intervals
$S_i$ of equal length, where $k$ is the distortion of either
$R^{\d}|_S$ or $R|_S$ and $m$ is subtorrentially big
(say $m<2^n$), the distortion obtained restricting to any interval $S_i$
will be bounded by $1+1/m$.

\end{rem}

\begin{rem} \label {remark 2}
 
We have the following estimate for the effect of the pullback of a subset of
$I_n$ by the central branch $R_n|_{I^0_n}$.
With total probability, for all $n$ sufficiently big,
if $X \subset I_n$ satisfies
$$
p_{\tg_n}(X|I_n)<\delta <n^{-b} \leq n^{-\b^{10\b}}
$$
then
$$
p_{\g_{n+1}}((R_n|_{I_{n+1}})^{-1}(X)|I_n)<\delta^{\a^3}<\delta^{-a}.
$$

Indeed, let $V$ be a $\delta^\a|I_{n+1}|$ neighborhood of $0$.  Then
$R_n|_{I_{n+1} \setminus V}$ has distortion bounded by $2 \delta^\a$.

Let $W \subset I_n$ an interval of size $\lambda |I_n|$.  Of course
$$
p_{\tg_n}(X \cap W|W)<\delta \lambda^{-\b}.
$$

Let's decompose each side of $I_{n+1} \setminus V$ as a union of
$n^\b \delta^{-\a}$ intervals.  Let $W$ be such an interval.
From Lemma \ref {away from the boundary},
it is clear that the image of $W$
covers at least $\delta^{2\a} n^{-3\b} |I_n|$.
It is clear then that
$$
p_{\tg_n}(X \cap R_n(W)|R_n(W))<\delta (\delta^{2\a}
n^{-3\b})^{-\b}=\delta^3 n^{-3 \b^2}.
$$
So we conclude (since the distortion of $R_n|_W$ is of order $1+n^{-3}$)
that
$$
p_{\g_{n+1}}((R_n|_{I_{n+1}})^{-1}(X) \cap W|W)<\delta^3 n^{-3 \b^2}
$$
(we use the fact that the composition of a $\g_{n+1}$-qs map with a map with
small distortion in $\tg_n$-qs).
Since
$$
p_{\g_{n+1}}(V|I_{n+1})<2 \delta^{\a^2},
$$
we get the required estimate.

\end{rem}

\subsection{More on trees} \label {more tree}

Let's see an application of the above remarks.

\begin{lemma}

With total probability, for all $n$ sufficiently big
$$
p_{\tg_n}((R^{\d}_n)^{-1}(X)|I^{\d}_n)<2^n p_{\g_n}(X|I_n).
$$

\end{lemma}

\begin{pf}

Decompose $I^{\d}_n$ in $n^{\ln(n)}$ intervals of equal length and
apply Remark \ref {remark 1}.
\end{pf}

By induction we get

\begin{lemma} \label {distortion product}

With total probability, for $n$ is big enough, if
$X_1,...,X_m \subset \Z \setminus \{0\}$
\begin{align*}
p_{\tg_n}(\d^{(n)}(x)=(j_1,...,j_m)&,j_i \in X_i|x \in I_n)\\
&\leq 2^{mn} \prod_{i=1}^m p_{\g_n}(j^{(n)}(x) \in X_i|x \in I_n).
\end{align*}

\end{lemma}

There is also a variation fixing the start of the sequence.

\begin{lemma}

With total probability, for $n$ is big enough, if
$X_1,...,X_m \subset \Z \setminus \{0\}$, and if $|\d|=(j_1,...,j_k)$
we have
\begin{align*}
p_{\tg_n}(\d^{(n)}(x)=(j_1,...,j_{m+k})&,j_{i+k} \in X_i|x \in I^\d_n)\\
&\leq 2^{mn} \prod_{i=1}^m p_{\g_n}(j^{(n)}(x) \in X_i|x \in I_n).
\end{align*}
In particular, with $\d=(\tau_n)$,
\begin{align*}
p_{\tg_n}(\d^{(n)}(x)=(\tau_n,j_1,...,j_m)&,j_i \in X_i|x \in
I^{\tau_n}_n)\\
&\leq 2^{mn} \prod_{i=1}^m p_{\g_n}(j^{(n)}(x) \in X_i|x \in I_n).
\end{align*}

\end{lemma}

The last part of the above Lemma will be often necessary
in order to apply PhPa1.

Sometimes we are more interested in the case where the $X_i$ are all equal.

Let $Q \subset \Z \setminus \{0\}$.  Let $Q(m,k)$ denote the set of
$\d=(j_1,...,j_m)$ such that $\#\{i|j_i \in Q\} \cap \{1,...,m\}
\geq k$.

Define $q_n(m,k)=
p_{\tg_n}(\cup_{\d \in Q(m,k)} I^{\d}_n|I_n)$.

Let $q_n=p_{\g_n}(\cup_{j \in Q} I^j_n|I_n)$.

\begin{lemma} \label {tree estimate}

With total probability, for $n$ large enough,
$$
q_n(m,k) \leq \binom {m} {k} (2^n q_n)^k.
$$

\end{lemma}

\begin{pf}

We have the following recursive estimates for $q_n(m,k)$:

$q_n(1,0)=1,q_n(1,1) \leq q_n \leq 2^n q_n$.

$q_n(m+1,0)=1,q_n(m+1,k+1) \leq q_n(m,k+1)+2^n q_n q_n(m,k)$.

Indeed, if $(j_1,...,j_{m+1}) \in Q(m+1,k+1)$ then either $(j_1,...,j_m) \in
Q(m,k+1)$ or $(j_1,...,j_m) \in Q(m,k)$ and $j_{m+1} \in Q$.
By the estimate $2^n$ on the distortion of all
non central branches, we get our result.
\end{pf}

We recall that by Stirling Formula,
$$
\binom {m} {q m}<\frac {m^{q m}} {(q m)!} <
\left (\frac {3} {q} \right )^{qm}.
$$

So we can get the following estimate.  For $q \geq q_n$,
$$
q_n(m,(6 \cdot 2^n) q m)<\left (\frac {1} {2}\right )^{(6 \cdot 2^n q m)}.
$$

It is also used in the following form.  If $q^{-1}>6 \cdot 2^n$ (it is
usually the case, since $q_n$ will be torrentially small)
$$
\sum_{k>q^{-2}} q_n(k,(6 \cdot 2^n) q k) <
2^{-n} q^{-1} \left (\frac {1} {2} \right )^{(6 \cdot 2^n) q^{-1}}.
$$

\section{Estimates on time}

Our aim in this section is to estimate the distribution of return times to
$I_n$:
they are concentrated around $c^{-1}_{n-1}$ up to an exponent in some range
given by $a$ and $b$.

The basic estimate is a large deviation estimate and is
proven in the next subsection
(Corollary \ref {large times estimate}) and states that for $k \geq 1$
the set of branches with time larger then
$k c_n^{-4\b}$ has capacity less then $e^{-k}$.

\subsection{A Large Deviation lemma for times}

Let $r_n(j)$ be such that $R_n|_{I^j_n}=f^{r_n(j)}$.  We will also use the
notation
$r_n(x)=r_n(j^{(n)}(x))$, the $n$-th return time of $x$ (there should be no
confusion for the reader, since we consistently use $j$ for an integer index
and $x$ for a point in the phase space.

Let
$$
A_n(k)=p_{\g_n}(r_n(x) \geq k|x \in I_n)
$$
Since $f$ restricted to the complement of $I_{n+1}$ is hyperbolic, from
Lemma \ref {hyperbol}, it is clear that
$A_n(k)$ decays exponentially with $k$.

Let $\zeta_n$ be the maximum $\zeta<c_{n-1}$ such that for all
$k>\zeta^{-1}$ we have
$$
A_n(k) \leq e^{-\zeta k}
$$
and finally let $\alpha_n=\min_{1 \leq m \leq n} \zeta_m$.

Our main result in this section is to estimate $\alpha_n$.  We will show
that with total probability, for $n$ big we have
$\alpha_{n+1} \geq c^{4\b}_n$.
For this we will have to do a simultaneous estimate for landing times,
which we define now.

Let $l_n(\d)$ be such that $L_n|_{I^\d_n}=f^{l_n(\d)}$.  We will also use
the notation $l_n(x)=l_n(\d^{(n)}(x))$.

Let
$$
B_n(k)=p_{\tg_n}(l_n(x)>k|x \in I_n).
$$

%Phase-estimate

\begin{lemma} \label {landing times}

If $k>c_n^{-2\b} \alpha_n^{-2\b}$ then
$$
B_n(k)<e^{-c_n^{2\b} \alpha_n^{2\b} k}.
$$

\end{lemma}

\begin{pf}

Let $k>c_n^{-2\b} \alpha_n^{-2\b}$ be fixed.  Let $m_0=\alpha_n^{2\b} k$.

Notice that by Lemma \ref {estimate on m}
$$
p_{\tg_n}(|\d^{(n)}(x)| \geq m_0|x \in I_n)\leq e^{-c_n^\b \alpha_n^{2\b} k}.
$$

Fix now $m<m_0$.  Let's estimate
$$
p_{\tg_n} (|\d^{(n)}(x)|=m,l_n(x)>k|x \in I_n).
$$

For each $\d=(j_1,...,j_m)$ we can associate a sequence of $m$
positive integers $r_i$ such that $r_i \leq r_n(j_i)$ and $\sum r_i=k$.
The average value of $r_i$ is at least $k/m$ so we conclude that
$$
\sum_{r_i \geq k/2m} r_i>k/2.
$$
Recall also that
$$
\frac {k} {2m}>\frac {1} {(2\alpha_n^{2\b})}>\alpha_n^{-1}.
$$

Given a sequence of $m$ positive integers $r_i$ as above we
can do the following estimate using Lemma \ref {distortion product}
\begin{align*}
p_{\tg_n} (\d^{(n)}(x)=(j_1,...,j_m)&,r_n(j_i)>r_i|x \in I_n)\\
&\leq 2^{mn} \prod_{j=1}^m p_{\g_n}(r_n(x) \geq r_j|x \in I_n)\\
&\leq
2^{mn} \prod_{r_j \geq \alpha_n^{-1}} p_{\g_n}(r_n(x) \geq r_j|x \in I_n)\\
&\leq 2^{mn} \prod_{r_j \geq k/2m} e^{-\alpha_n r_j}\\
&\leq 2^{mn} e^{-\alpha_n k/2}.
\end{align*}

The number of sequences of $m$ positive integers $r_i$ with sum $k$ is
\begin{align*}
\binom {k+m-1} {m-1} &\leq \frac {1}{(m-1)!}
(k+m-1)^{m-1}\\
&\leq \frac{1}{m!}(k+m)^m
\leq \left (\frac {2ek} {m} \right )^m.
\end{align*}

Notice that
\begin{align*}
2^{mn}\left (\frac {2ek} {m} \right)^m &\leq
\left (\frac {2^{n+3}k} {m} \right )^{\frac {m} {k2^{n+3}} k2^{n+3}}\\
&\leq \left (\frac {2^{n+3}k} {m_0} \right)^{\frac {m_0} {k2^{n+3}} k2^{n+3}}
(\text {since $x^{1/x}$ is decreases for $x>e$})\\
&\leq \left (\frac {2^{n+3}} {\alpha_n^{2\b}} \right )^{m_0} \leq
e^{\alpha_n^\b k}.
\end{align*}

So we can finally estimate
\begin{align*}
p_{\tg_n} (|\d^{(n)}(x)|=m,l_n(x) \geq k|x \in I_n) &\leq 2^{mn}
\left ( \frac {2ek} {m} \right )^m
e^{-\alpha_n k/2}\\
&< e^{(\alpha_n^{\b-1}-1/2)\alpha_n k}.
\end{align*}

Summing up on $m$ we get
\begin{align*}
p_{\tg_n} (|\d^{(n)}(x)|<m_0&,l_n(x) \geq k|x \in I_n)\\
&\leq m_0 e^{(\alpha_n^{\b-1}-1/2)\alpha_n k}\\
&<e^{(2 \alpha_n^{\b-1}-1/2) \alpha_n k}
(\text {since}
\frac {\ln(m_0)} {k} \leq \frac {\ln(k)} {k} \leq
\alpha_n^\b)\\
&\leq e^{-\alpha_n k/3}.
\end{align*}

As a direct consequence we get 
$$
B_n(k)<e^{-\alpha_n k/3}+
e^{-c_n^\b \alpha_n^{2\b} k}<e^{-c_n^{2\b} \alpha_n^{2\b} k}.
$$
\end{pf}

%Phase-parameter

\begin{lemma} \label {estimate on v_n}

With total probability, for $n$ large enough,
$$
v_{n+1}<c^{-3\b}_n \alpha_n^{-3\b}.
$$

\end{lemma}

\begin{pf}

Let $\d$ such that $R_n(0) \in C^{\d}_n$.

Using Lemma \ref {landing times}, with total probability,
for $n$ large enough,
$l_n(\d)<n \alpha_n^{-2\b} c_n^{-2\b}$, and
$v_{n+1}<v_n+n \alpha_n^{-2\b} c_n^{-2\b}$.
Using torrential (and monotonic)
decay of $\alpha_n c_n$ we get for $n$ large
enough $v_{n+1}<c_n^{-3\b} \alpha_n^{-3\b}/2$.
\end{pf}

%Phase-estimate

\begin{lemma}

With total probability, for $n$ large enough,
$$
\alpha_{n+1} \geq \min \{\alpha_n^{4\b}, c_n^{4\b}\}.
$$

\end{lemma}

\begin{pf}

Let $k \geq \max\{\alpha_n^{-4\b},c_n^{-4\b}\}$.
From Lemma \ref {estimate on v_n} one
immediately sees that if $r_{n+1}(j) \geq k$ then $R_n(I^j_{n+1}))$ is
contained on some $C^{\d}_n$ with $l_n(\d) \geq
k/2 \geq \alpha_n^{-2\b} c_n^{-2\b}$.

Applying Lemma \ref {landing times} we have
$B_n(k/2)<e^{-\alpha_n^{2\b} c_n^{2\b} k/2}$.

Applying Remark \ref {remark 2} we get
$$
A_{n+1}(k)<e^{-k \alpha_n^{2\b}
c_n^{2\b}/\a^3}<e^{-k \min\{\alpha_n^{4\b},c_n^{4\b}\}}.
$$
\end{pf}

Since $c_n$ decreases torrentially, we get

%Phase-estimate

\begin{cor} \label {large times estimate}

With total probability, for $n$ large enough $\alpha_{n+1} \geq c_n^{4\b}$.

\end{cor}

\begin{rem} \label {v_n}

In particular, using Lemma \ref {estimate on v_n}, for $n$ big,
$v_n<c_{n-1}^{-6\b}$.

\end{rem}

\subsection{Consequences}

The Lemma bellow is just a convenient way to summarize our results on the
distribution of times.  We also take the opportunity to state them in terms
of constants $a$ and $b$.

\begin{lemma} \label {large times}

With total probability, for all $n$ sufficiently large we have

\begin{enumerate}

\item $p_{\tg_n}(l_n(x)<c_n^{-s}|x \in I_n)<c_n^{\sqrt \a-s} <
c_n^{a-s}$, with $s>0$,

\item $p_{\tg_n}(l_n(x)<c_n^{-s}|x \in I^{\tau_n}_n)<c_n^{a-s}$, with $s>0$,  

\item $p_{\tg_n}(l_n(x)>c_n^{-s}|x \in I_n)<e^{-c_n^{b-s}}$, with $s>b$,

\item $p_{\tg_n}(l_n(x)>c_n^{-s}|x \in I^{\tau_n}_n)<e^{-c_n^{b-s}}$,
with $s>b$,

\item $p_{\g_n}(r_n(x)<c_{n-1}^{-s}|x \in I_n)<c_{n-1}^{\sqrt \a-s} <
c_{n-1}^{a-s}$, with $s>0$,

\item $p_{\g_n}(r_n(x)>c_{n-1}^{-s}|x \in I_n)<e^{-c_{n-1}^{\sqrt b-s}} <
e^{-c_{n-1}^{b-s}}$ with $s>b$.

\item $c_{n-1}^{-a}<r_n(\tau_n)<c_{n-1}^{-b}$.

\item $c_{n-1}^{-a}<v_n<c_{n-1}^{-b}$.

\end{enumerate}

\end{lemma}

\begin{pf}

The first and third estimates are contained in
Lemma \ref {estimate on m} (after noticing $l_n(x) \geq |\d^{(n)}(x)|$)
and the second follows from Lemma \ref {landing times}.  $6$ is contained in
Corollary \ref {large times estimate}, while $5$ can be
obtained by pulling back estimate $1$ (using Remark \ref {remark 2}).
$5$ and $6$ imply $7$ by PhPa2.  In view of $7$, $4$ follows from the
proof of Lemma \ref {landing times} ($7$ is needed to avoid
$r_n(\tau_n)$ to be too big).  $8$ is easily obtained from
Lemma \ref {growth of s_n} and Remark \ref {v_n}.
\end{pf}

\begin{rem} \label {precise rate of c_n}

It is clear (using Lemma \ref {away from the boundary})
that $8$ in Lemma \ref {large times} implies that with total
probability, for all $n$ big enough, $c_{n+1} \geq
4^{-c^{-b}_n}$.  This implies together with Corollary \ref {c_n
torrential} that
$$
a<\frac {\ln(\ln(c_n^{-1}))} {\ln(c_{n-1}^{-1})}<b,
$$
so $c^{-1}_n$ grows torrentially.

\end{rem}

\section{Dealing with hyperbolicity}

In this section we show by an inductive proccess that the great
majority of branches are reasonably hyperbolic.
In order to do that, in the
following subsection, we define some classes of branches with
`good' distribution of times and which are not too close to the critical
point.  The definition of `good' distribution of
times has an inductive
component: they are composition of many `good' branches of the previous
level.  The fact that most branches are good is related to the validity of
some kind of Law of Large Numbers estimate.

\subsection{Some kinds of branches and landings}

\subsubsection{Standard landings}

We define the set of standard landings at time $n$,
$LS(n) \subset \Omega$ as the set
of all $\d=(j_1,...,j_m)$ satisfying the following.

\begin{description}

\item[LS1] ($m$ not too small or large) $c^{-a/2}_n<m<c^{-2b}_n$,

\item[LS2] (No very large times) $r_n(j_i)<c^{-3b}_{n-1}$ for all $i$.

\item[LS3] (Short times are sparse in not too small initial segments)
For $c^{-2b}_{n-1} \leq k \leq m$
$$
\#\{r_n(j_i)<c^{-a/2}_{n-1}\} \cap \{1,...,k\}<
(6 \cdot 2^n) c^{a/2}_{n-1} k.
$$
  
  \end{description}

We also define the set of fast landings at time $n$,
$LF(n) \subset \Omega$  by the following conditions
\begin{description}

\item[LF1] ($m$ small) $m<c^{-a/2}_n$.

\item[LS2] (no very large times) $r_n(j_i)<c^{-3b}_{n-1}$ for all $i$.
  
  \end{description}

%Phase estimate

\begin{lemma} \label {standard landing}

With total probability, for all $n$ sufficiently big,

\begin{enumerate}

\item $p_{\tg_n}(\d^{(n)}(x) \notin LS(n)|x \in I_n)<c_n^{a/3}/2$,

\item $p_{\tg_n}(\d^{(n)}(x) \notin
LS(n) \cup LF(n)|x \in I_n)<c_n^{n^2}/2$,

\item $p_{\tg_n}(\d^{(n)}(x) \notin LS(n)|x \in I^{\tau_n}_n)<c_n^{a/3}/2$,

\item $p_{\tg_n}(\d^{(n)}(x) \notin LS(n) \cup LF(n)|x \in
I^{\tau_n}_n)<c_n^{n^2}/2$.

\end{enumerate}
  
\end{lemma}

\begin{pf}

The proof is immediate from our time estimates, which can be applyed
directly to estimate the losses of LS1, LS2, LF1, or under large
deviation form for LS3, following \S \ref {more tree}.
Let's estimate the complement of the sets which do not satisfy some
properties:

\begin{description}

\item [LS1] $c_{n-1}^{2a/5}$,

\item [LS1+LF1] $e^{c_n^b} \ll c_n^{n^2}$,

\item [LS2] $e^{c_{n-1}^{-2b}} \ll c_n^{n^2}$,

\item [LS3] let $q=(6 \cdot 2^n) c_{n-1}^{a/2}$, we get
$$
\sum_{k>c_{n-1}^{-2b}} \left (\frac {1} {2} \right )^{qk} \ll c_n^{n^2}.
$$

\end{description}

This gives immediately $1$ and $2$.  For $3$ and $4$ the same estimates
hold, using (for LS2) the estimate on $r_n(\tau_n)$.
\end{pf}

\subsubsection{Very good returns, bad returns and excellent landings}

Define the set of very good returns,
$VG(n_0,n) \subset \Z \setminus \{0\}$, $n_0,n \in \N$ and
the set of bad returns,
$B(n_0,n) \subset \Z \setminus \{0\}$, $n_0,n \in \N$,
$n \geq n_0$ by induction as follows.  We let
$VG(n_0,n_0)=\Z \setminus \{0\}$,$B(n_0,n_0)= \emptyset$ and supposing
$VG(n_0,n)$ and $B(n_0,n)$ defined, define the set of excellent landings
$LE(n_0,n) \subset LS(n)$ satisfying the following extra assumptions.   

\begin{description}

\item[LE1] (Not very good moments are sparse in not too small initial
segments)
For all $c^{-2b}_{n-1}<k \leq m$
$$
\#\{i|j_i \notin VG(n_0,n)\} \cap \{1,...,k\}<
6 \cdot 2^n c^{a^2}_{n-1} k,
$$
\item[LE2] (Bad moments are sparse in not too small initial segments)
For all $c^{-1/n}_{n}<k \leq m$
$$
\#\{i|j_i \notin B(n_0,n) \} \cap \{1,...,k\}<
6 \cdot 2^n c^n_{n-1} k,
$$
  
\end{description}

We define $VG(n_0,n+1)$ as the set of $j$ such that
$R_n(I^j_{n+1})=C^{\d}_n$ with $\d \in LE(n_0,n)$ and
the extra condition.

\begin{description}

\item[VG] (distant from $0$) The distance of $I^j_{n+1}$ to $0$ is bigger
than $c_n^{n^2}|I_{n+1}|$.

\end{description}

And we define $B(n_0,n+1)$ as the set of $j \notin VG(n_0,n+1)$ such that
$R_n(I^j_{n+1})=C^{\d}_n$ with $\d \notin LF(n)$.

%Phase estimate

\begin{lemma} \label {induction step}
With total probability, for all $n_0$ sufficiently big,
\begin{enumerate}
\item $p_{\g_n}(j^{(n)}(x) \notin VG(n_0,n)|x \in I_n)<c_{n-1}^{a^2},$
\item $p_{\g_n}(j^{(n)}(x) \in B(n_0,n)|x \in I_n)<2 c_{n-1}^{2 n}$,
\item $p_{\tg_n}(\d^{(n)}(x) \notin LE(n_0,n)|x \in I_n)<c_n^{2a/5}$,
\item $p_{\tg_n}(\d^{(n)}(x) \notin LE(n_0,n) \cup LF(n)|x \in I_n) <
c_n^{bn}$,
\item $p_{\tg_n}(\d^{(n)}(x) \notin LE(n_0,n)|x \in I^{\tau_n}_n) <
c_n^{2a/5}$.
\end{enumerate}

\end{lemma}

\begin{pf}

First notice that the validity for a given value of $n$ of $1$ and $2$
implies $3$, $4$ and $5$, using the large deviation
technique of \S \ref {more tree}.
More precisely, we estimate the complement

\begin{description}

\item [LE1] let $q=6 \cdot 2^n c_{n-1}^{a^2}$,
$$
\sum_{k>c_{n-1}^{-2b}} \left (\frac {1} {2} \right )^{qk} \ll c_n^{n^2},
$$

\item [LE2] let $q=6 \cdot 2^n c_{n-1}^n$,
$$
\sum_{k>c_n^{-1/n}} \left (\frac {1} {2} \right )^{qk} \ll c_n^{n^2}.
$$

\end{description}

The validity of $3$ and $4$ imply
$1$ and $2$ for $n+1$ by pulling back, using Remark \ref {remark 2}.
The proof then follows by induction (the first step is trivial).
\end{pf}

Using PhPa2 we get

%Phase-parameter

\begin{lemma} \label {lands very good}

With total probability, for all
$n_0$ big enough, for all $n$ big enough,
$\tau_n \in VG(n_0,n)$.

\end{lemma}

We also have the following trivial estimate for a very good return time

%Phase-estimate

\begin{lemma} \label {very good return time}

With total probability, for all
$n_0$ big enough and for all $n \geq n_0$,
if $j \in VG(n_0,n+1)$ then
$$
m <r_{n+1}(j)<m c^{-4b}_{n-1},
$$
where as usual, $m$ is such that $R_n(I^j_{n+1})=C^{\d}_n$ and
$\d=(j_1,...,j_m)$.

\end{lemma}

\begin{pf}

The estimate from below is obvious, the estimate from above follows from LS2
and the estimate on $v_n$.
\end{pf}

And the following estimate for returns which are not very good or bad

\begin{lemma} \label {not very good or bad}

With total probability for all $n_0$ sufficiently big, if $n>n_0$,
if $j \notin VG(n_0,n) \cup B(n_0,n)$ then
$r_n(j)<c_{n-1}^{-a/2} c_{n-2}^{-4b}$.

\end{lemma}

\begin{pf}

If $j \notin VG(n_0,n) \cup B(n_0,n)$ then $R_{n-1}(I^j_n) \in LF(n_0,n-1)$. 
The estimate follows since a branch in $LF(n_0,n-1)$ has time bounded by
$c_{n-1}^{a/2}c_{n-2}^{-3b}$ and $v_{n-1}<c_{n-2}^{-b}$.
\end{pf}

Let $j \in VG(n_0,n+1)$.  We can write
$R_{n+1}|_{I^j_{n+1}}=f^{r_{n+1}(j)}$, that is, a big iterate of $f$.  One
may consider which proportion of those iterates belong to very good
branches of the previous level.  More generally, we can truncate the return
$R_{n+1}$, that is, we may consider
$k<r_{n+1}(j)$ and ask which proportion of iterates up to $k$ belong to very
good branches.

%Phase-estimate

\begin{lemma} \label {very good partial time}

With total probability, for all
$n_0$ big enough and for all $n \geq n_0$,
the following holds.

Let $j \in VG(n_0,n+1)$, as usual let
$R_n(I^j_{n+1})=I^{\d}_n$ and
$\d=(j_1,...,j_m)$.  Let $m_k$ be biggest possible with
$$
v_n+\sum_{j=1}^{m_k} r_n(j_i) \leq k
$$
(the amount of full returns to level $n$ before time $k$) and let
$$
\beta_k=\sum_{\ntop {1 \leq i \leq m_k,}
{j_i \in VG(n_0,n)}} r_n(j_i).
$$
(the total time spent in full returns to level $n$
which are very good before time $k$)
Then $1-\beta_k/k<c_{n-1}^{a^2/3}$ if $k>c_n^{-2/n}$.

\end{lemma}

\begin{pf}

Let's estimate first the time $i_k$ which is not spent on full returns:
$$
i_k=k-\sum_{j=1}^{m_k} r_n(j_i).
$$
This corresponds exactly to $v_n$ plus some incomplete part of the return
$j_{m_{k+1}}$.  This part can be bounded by $c_{n-1}^{-b}+c_{n-1}^{-3b}$ 
(use the estimate of $v_n$ and LS2 to estimate
the incomplete part).

Using LS2 we conclude now that
$$
m_k>(k-c_{n-1}^{-b})c_{n-1}^{3b}>c_n^{-1/n}
$$
so $m_k$ is not too small.

Let's now estimate the contribution $h_k$ from bad full returns $j_i$.  
The number of such returns
must be less than $c_{n-1}^{n/2} m_k$, using LS2 their total time is at most
$c_{n-1}^{(n/2)-3b} m_k<m_k$.

The non very good full returns on the other hand can be estimated by LE1
(given the estimate on $m_k$), they are at most $c_{n-1}^{a^2} m_k$.
So we can estimate the total time $l_k$ of non very good full returns   
with time less then $c_{n-1}^{-a/2} c_{n-2}^{-4b}$ by
$$
c_{n-1}^{a^2}c_{n-1}^{-a/2} c_{n-2}^{-4b} m_k,
$$
while $\beta_k$ can be estimated from below by
$$
(1-c^{a/4}_{n-1}) c^{-a/2}_{n-1} m_k.
$$      

It is easy to see then that $i_k/\beta_k \ll c^{a/5}_{n-1}$,
$h_k/\beta_k \ll c^{a/5}_{n-1}$.
We also have
$$
l_k/\beta_k<2 c_{n-1}^{a^2/2}.
$$
So $(i_k+h_k+l_k)/\beta_k$ is less then $c^{a^2/3}_{n-1}$.
Since $i_k+h_k+l_k+\beta_k=k$ we have $1-\beta_k/k<(i_k+h_k+l_k)/\beta_k$.
\end{pf}

\subsubsection{Cool landings}

Let's define the set of cool landings
$LC(n_0,n) \subset \Omega$, $n_0,n \in \N$,
$n \geq n_0$ as the set of all $\d=(j_1,...,j_m)$ in $LE(n_0,n)$
satisfying.

\begin{description}

\item[LC1] (Starts very good) $j_i \in VG(n_0,n)$, $1 \leq i \leq
c_{n-1}^{-a^2/2}$.

\item[LC2] (Short times are sparse in not too small initial segments)
For $c^{-a/2}_{n-1} \leq k \leq m$
$$
\#\{r_n(j_i)<c^{-a/2}_{n-1}\} \cap \{1,...,k\}<
(6 \cdot 2^n) c^{a/3}_{n-1} k,
$$
  
\item[LC3] (Not very good moments are sparse in not too small initial
segments)
For all $c^{-a^2/4}_{n-1}<k \leq m$
$$
\#\{i|j_i \notin VG(n_0,n)\} \cap \{1,...,k\}<
(6 \cdot 2^n) c^{a^2}_{n-1} k,
$$
  
\item[LC4] (Bad times are sparse in not too small initial segments)
For $c_{n-1}^{-n/3} \leq k \leq m$
$$
\#\{i|j_i \in B(n_0,n)\} \cap \{1,...,k\}<
(6 \cdot 2^n) c^{n/6}_{n-1} k,
$$
  
\item[LC5] (Starts with no bad times) $j_i \notin B(n_0,n)$,
$1 \leq i \leq c_{n-1}^{-n/2}$.

\end{description}

%Phase-estimate

As usual we obtain

\begin{lemma}

With total probability, for all
$n_0$ sufficiently big and all $n \geq n_0$,
$$
p_{\tg_n}(\d^{(n)}(x) \notin LC(n_0,n)|x \in I_n)<c_{n-1}^{a^2/3}
$$
and for all $n$ big enough
$$
p_{\tg_n}(\d^{(n)}(x) \notin LC(n_0,n)|x \in I^{\tau_n}_n)<c_{n-1}^{a^2/3}.
$$

\end{lemma}

\begin{pf}

\begin{description}

\item [LC1] $2^n c_{n-1}^{a^2/2} c_{n-1}^{a^2}<c_{n-1}^{a^2/3}$,

\item [LC2] let $q=6 \cdot 2^n c_{n-1}^{a/3}$,
$$
\sum_{k>c_n^{-a/2}} \left (\frac {1} {2} \right )^{qk} \ll c_{n-1}^{n^2}.
$$

\item [LC3] let $q=6 \cdot 2^n c_{n-1}^{a^2}$,
$$
\sum_{k>c_{n-1}^{-a^2/4}} \left (\frac {1} {2} \right )^{qk} \ll
c_{n-1}^{n^2},
$$

\item [LC4] let $q=6 \cdot 2^n c_{n-1}^{n/6}$,
$$
\sum_{k>c_n^{-n/3}} \left (\frac {1} {2} \right )^{qk} \ll c_n^{n^2}.
$$

\item [LC5] $2^n c_{n-1}^{-n/2} c_{n-1}^{2 n} \ll c_{n-1}^n$.

\end{description}

This gives the first estimate, for the second we must use Lemma \ref {lands
very good} and the estimates on $r_n(\tau_n)$.

\end{pf}

Using PhPa1 we get

%Phase-parameter

\begin{lemma} \label {lands cool}

With total probability, for all
$n_0$ big enough, for all $n$ big enough we
have $R_n(0) \in LC(n_0,n)$.

\end{lemma}

\subsection{Hyperbolicity}

\subsubsection{Preliminaries}

For $j \neq 0$, we define
$$
\lambda^{(j)}_n=\inf_{x \in I^j_n} \frac
{\ln(|R_n'(x)|)} {r_n(j)}.
$$
And $\lambda_n=\inf_{j \neq 0} \lambda^{(j)}_n$.  As a consequence of the
exponential estimate on distortion, together with hyperbolicity of $f$ in
the complement of $I^0_n$ we immediately have the following.

%Phase-estimate

\begin{lemma}

With total probability, for all $n$ sufficiently big, $\lambda_n>0$.

\end{lemma}

\subsubsection{Good branches}

We define the set of good returns
$G(n_0,n) \subset \Z \setminus \{0\}$, $n_0,n \in \N$,
$n \geq n_0$ as the set of all $j$ such that

\begin{description}

\item[G1] (hyperbolic return)
$$
\lambda^{(j)}_n \geq \lambda_{n_0} \frac {1+2^{n_0-n}} {2},
$$

\item[G2] (hyperbolicity in partial return)
for $c_{n-1}^{-3/(n-1)} \leq k \leq r_n(j)$ we have
$$
\inf_{I^j_n} \frac {\ln(|Df^k|)} {k} \geq
\lambda_{n_0} \frac {1+2^{n_0-n+1/2}} {2}-c_n^{2/(n-1)}.
$$

\end{description}

Notice that since $c_n$ decreases torrentially, for $n$ sufficiently big G2
implies for $c_{n-1}^{-3/(n-1)} \leq k \leq r_n(j)$ we have
$$
\inf_{I^j_n}\frac {\ln(|Df^k|)} {k} \geq
\lambda_{n_0} \frac {1+2^{n_0-n}} {2}.
$$

%Phase-parameter

\begin{lemma} \label {very good is good}

With total probability, for $n_0$ big enough, for all $n>n_0$,
$VG(n_0,n) \subset G(n_0,n)$.

\end{lemma}

\begin{pf}

Let's prove that if G1 is satisfied
for all $j \in VG(n_0,n)$, then $VG(n_0,n+1) \subset G(n_0,n+1)$.  Let's fix
such a $j$.  Notice
that by definition of $\lambda_{n_0}$ the hypothesis is satisfied for $n_0$.

Let $a_k=\inf_{I^j_n} \ln(|Df^k|)/k$.

Recall the estimate from Lemma \ref {very good
return time} for $t=r_{n+1}(j)$, we as usually let
$R_n(I^j_{n+1})=C^{\d}_n$, $\d=(j_1,...,j_m)$.
Let's say that $j_i$ was completed before $k$ if
$v_n+r_n(j_1)+...+r_n(j_i) \leq k$.  We let the queue be defined as
$$
q_k=\inf_{C^{\d}_n} \ln (|Df^{k-r} \circ f^r|)
$$
where
$r=v_n+r_n(j_1)+...+r_n(j_{m_k})$ with $j_{m_k}$ the last complete return.

Recall that $v_n < c_{n-1}^{-b}$ and from VG it is clear
$|R'_n|_{I^j_{n+1}}|>c_n^{n^2}$.  Notice also that using Lemma \ref
{distortion}, for any $k_0 \leq m$,
the derivative of $R^{k_0}_n$ in $C^\d_n$ is at least $2^{k_0}$.  So for
$m_0=c_{n-1}^{-2b}$ we have that the derivative of $R_n^{k_0+1}$ in
$I^j_{n+1}$ is at least $1$.  Notice that from LS2
$$
k_0=\sum_{i=1}^{m_0} r_n(j_i)<c_{n-1}^{-2b} c_{n-1}^{-3b} \ll k.
$$

It is clear that any complete return
before $k$ gives derivative at least $1$ from Lemma \ref {distortion}.
The queue can
be bound by $\ln(c_n c_{n-1}^3)$ using Lemma \ref {lower bound}.  We have
$-q_k/k \ll c_n^{2/n}$.

Now we use Lemma \ref {very good partial time} and get
\begin{align*}
a_k &> \frac {\beta_k-k_0} {k} \frac {\lambda_{n_0}(1+2^{n_0-n})} {2}-
\frac {-q_k} {k}\\
&\geq \frac {\lambda_{n_0}(1+2^{n_0-n-1/2})} {2}-\frac {-q_k} {k}.
\end{align*}

Which gives G2.  If $k=r_{n+1}(j)$, $q_k=0$ which gives G1.
\end{pf}

\subsubsection{Hyperbolicity in cool landings}

%Phase-estimate

\begin{lemma} \label {cool hyperbolicity}

With total probability, if $n_0$ is sufficiently big, for all $n$
sufficiently big, if $\d \in LC(n_0,n+1)$ then for all
$c_{n-1}^{-4/(n-1)}<k \leq l_n(\d)$,
$$
\inf_{C^{\d}_n} \frac {\ln(|Df^k|)} {k} \geq
\frac {\lambda_{n_0}} {2}.
$$

\end{lemma}

\begin{pf}

Fix such $\d \in LC(n_0,n+1)$.

Let
$$
a_k=\inf_{C^{\d}_n} \frac {\ln(|Df^k|)} {k}.
$$

Consider the sequence $r_i=r_n(j_i)$ where as usual
$\d=(j_1,...,j_m)$.

As in Lemma \ref {very good partial time}, we define $m_k$ as the biggest
such that
$$
\sum_{i=1}^{m_k} r_i \leq k.
$$
We define
$$
\beta_k=\sum_{\ntop {1 \leq i \leq m_k,}
{j_i \in VG(n_0,n+1)}} r_i,
$$
(counting the time up to $k$ spent in complete very good returns)
and
$$
i_k=k-\sum_{i=1}^{m_k} r_i.
$$
(counting the time in the incomplete return at $k$).

Let's then consider two cases: small $m_k$ ($m_k<c_{n-1}^{-a^2/2}$)
and otherwise.

The idea of the first case is that all full returns are very good by LC1,
and the incomplete time is also part of a very good return.

Since full very good returns are very hyperbolic by G1 and
very good returns are good, we
just have to worry with possibly losing hyperbolicity in the incomplete
time.  To control this, we introduce the queue
$q_k=\inf_{C^{\d}_{n-1}} \ln(|Df^{i_k} \circ f^{k-i_k}|)$.
We have $-q_k<-\ln(c_{n-1}^{1/3}
c_{n-1}^3)$ by Lemma \ref {lower bound}

If the incomplete time is big (more than $c_{n-1}^{-4/(n-1)}$), we can use
G2 to estimate the hyperbolicity of the incomplete time (which is part of a
very good return).  The reader can easily check the estimate in this case.

If the incomplete time is not big, we can not use G2 to estimate $q_k$, but
in this case $i_k$ is much less than $k$: since $k>c_{n-1}^{-4/(n-1)}$, at
least one return was completed ($m_k \geq 1$), and since it must be very
good we conclude that $k>c_{n-1}^{-a/2}$ by LS1, so
$$
a_k>\lambda_{n_0} \frac {(1+2^{n-n_0})} {2} \cdot \frac {k-i_k} {k}-
\frac {-q_k} {k}>\frac {\lambda_{n_0}} {2}.
$$

Let's consider now the case $m_k>c_{n-1}^{-a^2/2}$.
For an incomplete time we still have $-q_k<-\ln(c_n c_{n-1}^3)$,
so $-q_k/k<c_{n-1}^{-a^2/3}$.

Arguing as in Lemma \ref {very good partial time}, we split $k-\beta_k-i_k$
(time of full returns which are not very good) in part relative to
bad returns $h_k$ and in part relative to fast returns (not very good or
bad) $l_k$.  Using LC4 and LC5
we get
$$
h_k<c_{n-1}^{-3b} c_{n-1}^{n/7} m_k,
$$
and using LC1 and LC3 we have
$$
l_k<c_{n-1}^{-a/2} c_{n-2}^{-4b} (6 \cdot 2^n) c_{n-1}^{a^2} m_k,
$$
By LC1 and LC3 we have $\beta_k>c_{n-1}^{-a/2} m_k/2$ so we have
$$
\frac {h_k+l_k} {\beta_k}<c_{n-1}^{a^2/2}.
$$

By LC1, $\beta_k>c_{n-1}^{a/2} m_k$.
Now $\beta_k>c_{n-1}^{-a/2}c_{n-1}^{-a^2}/2$ by hypothesis on $m_k$.
If $i_k<c_{n-1}^{-a/2}c_{n-2}^{-4b}$,
$i_k/\beta_k$ is small and we are done.  Otherwise by Lemma \ref {not very
good or bad}, $i_k$ must be either very good or bad.  If $i_k$ is very
good we can reason
as before that $G2$ applies for the estimate of the queue and we are done. 

If $i_k$ is bad, by LC5 we have that $m_k>c_{n-1}^{n/2}$,
but $i_k<c_{n-1}^{3b}$ by LS2, so $i_k/\beta_k$ is
very small again and we are done.
\end{pf}

\section{Main Theorems}

\subsection{Proof of Theorem A}

%Phase-estimate

We must show that with total probability, $f$ is Collet-Eckmann.  We will
use the estimates on hyperbolicity of cool landings to show that if the
critical point always fall in a cool landing then there is uniform control
of the hyperbolicity along the critical orbit.

Let
$$
a_k=\frac {\ln(|Df^k(f(0))|)} {k}
$$
and $e_n=a_{v_n-1}$.

It is easy to see that if $n_0$ is big enough such that both Lemmas \ref
{lands cool} and \ref {very good is good} are valid
and using Lemma \ref {away from the boundary} so that
$|R_n(0)|>|I_n|/2^n$ for $n$ big
enough, we obtain using Lemma \ref {cool hyperbolicity} that
$$
e_{n+1} \geq e_n \frac {v_n-1} {v_{n+1}-1}+\frac {\lambda_{n_0}} {2} \frac
{v_{n+1}-v_n} {v_{n+1}-1}
$$
and so
$$
\lim \inf e_n \geq \frac {\lambda_{n_0}} {2}.
$$

Let now $v_n-1<k<v_{n+1}-1$.  Define
$$
q_k=\ln(|Df^{k-v_n}(f^{v_n}(0))|).
$$

If $k<v_n+c_{n-1}^{-4/(n+1)}$ from LC1 and LS1 we know that the
time of the $R_n$ branch of $R_n(0)$ is at least $c_{n-1}^{-a/2}$,
so $k$ is in the middle of this branch.  Using  $|R_n(0)|>|I_n|/2^n$, we
get that $-q_k<-\ln(2^{-n} c_{n-1} c_{n-1}^{-3})$.
We then get from $v_n>c_{n-1}^{-a}$ that
$$
a_k \geq e_n \frac {v_n-1} {k}-\frac {-q_k} {k}>(1-1/2^n)e_n-1/2^n.
$$

If $k>v_n+c_{n-1}^{-4/(n+1)}$ using Lemma \ref {cool hyperbolicity}
we get
$$
a_k \geq e_n \frac {v_n-1} {k}+\frac {\lambda_{n_0}} {2} \cdot
\frac {k-v_n+1} {k}.
$$

Those two estimates imply that $\lim \inf a_k \geq \lambda_{n_0}/2$ and so
$f$ is Collet-Eckmann.

\subsection{Proof of Theorem B}

We must obtain, with total probability, upper and lower (polynomial)
bounds for the recurrence of the critical orbit.  It will be easier to
first study the recurrence with respect to iterates of return branches,
and then estimate the total time of those iterates.

\subsubsection{Recurrence in terms of return branches}

\begin{lemma}

With total probability, for $n$ big enough and for
$1 \leq i \leq c^{-2}_{n-1}$,
$$
\frac {\ln(|R_n^i(0)|)} {\ln(c_{n-1})} < b^2.
$$

\end{lemma}

\begin{pf}

Notice that due to torrential (and monotonic) decay of $c_n$, we can
estimate $|I_n|=c_{n-1}^{1+\delta_n}$, with $\delta_n$ decaying torrentially
fast.

For $i=1$ it follows from Lemma \ref {away from the boundary}.

Let $X \subset I_n$
be a $c_{n-1}^{b^2}$ neighborhood of $0$.
For $n$ big, we can estimate (due to the relation between $|I_n|$ and
$c_{n-1}$)
$$
\frac {|X|} {|I_n|}<c_{n-1}^{b^2-2}
$$
(we of course consider $X_n$ a union of $C^\d_n$, so that its
size is near the required size, the precision is high enough for our
purposes due to Lemma \ref {approx}).

We have to make sure that the critical point does not land in $X$ for
$1<i \leq c_{n-1}^{-2}$.  This requirement can
be translated on $R_n(0)$ not belonging to a certain set
$Y \subset I_n$ such that
$$
Y=\bigcup_{1 \leq |\d|<c_{n-1}^{-2}} (R_n^{\d})^{-1}(X).
$$

It is clear that
$$
p_\g(I^{\tau_n}_n \cap Y|I^{\tau_n}_n) \leq
c_{n-1}^{-2}c_{n-1}^{(b^2-2)/b}<
c_{n-1}^{b-3}.
$$

Applying PhPa1, the probability that for
$1 \leq i \leq c^{-2}_{n-1}$ we have $|R^i_n(0)|<c^{b^2}_{n-1}$
is at most $c_{n-1}^{b-3}$, which is summable.
\end{pf}

\begin{lemma}

With total probability, for $n$ big enough and for
$c_{n-1}^{-2}<i \leq s_n$,
$$
\frac {\ln(|R_n^i(0)|} {\ln(c_{n-1}^{-1})}<
b^2 \left( 1+\frac {\ln(i)} {\ln(c_{n-1}^{-1})} \right).
$$

\end{lemma}

\begin{pf}

The argument is the same as for the previous Lemma, but the decomposition
has a slight different geometry.

For $j \geq 0$, let $X_j \subset I_n$
be a
$$
c_{n-1}^{b^2 (b^{j+1})}
$$
neighborhood of $0$ (approximated as union of $C^\d_n$, if $X_j \subset
I_{n+1}$, it won't be relevant for the proof, the reader can take then
$X_j=\emptyset$).  Let $Y_j \subset I_n$ be such that
$$
Y_j=\bigcup_{c_{n-1}^{-b^j} \leq |\d|<
c_{n-1}^{-b^{j+1}}} (R_n^{\d})^{-1}(X_j).
$$

It is clear that
$$
p_\g(I^{\tau_n}_n \cap Y_j|I^{\tau_n}_n) \leq
c_{n-1}^{-b^{j+1}}c_{n-1}^{b^{j+2}}<
c_{n-1}^{b^{j+1}}
$$
and
$$
p_\g(I^{\tau_n}_n \cap \cup
Y_j|I^{\tau_n}_n)<\sum_{j=0}^\infty
c_{n-1}^{-b^j}<2 c_{n-1}.
$$

Applying PhPa1, with total probability,
the critical point does not belong to any $Y_j$. 
Let's see what this means.  If $c_{n-1}^{-b^j}<i \leq c_{n-1}^{-b^{j+1}}$,
$R_n^i(0) \notin X_j$, that is $|R_n^i(0)|>c_{n-1}^{b^{j+3}}$.  This clearly
implies the statement.
\end{pf}

\begin{cor} \label {distance}

With total probability, for $n$ big enough and for
$1 \leq i \leq s_n$,
$$
\frac {\ln(|R_n^i(0)|)} {\ln(c_{n-1})} <
b^2 \left( 1+\frac {\ln(i)} {\ln(c_{n-1}^{-1})} \right).
$$

\end{cor}

\subsubsection{Total time of full returns}

For $1 \leq i \leq s_n$, let $k_i$ such that $R_n^i|_{I_{n+2}}=f^{k_i}$.

\begin{lemma}

With total probability, for $n$ big enough and for
$c_{n-1}^{-1}<i<s_n$, $k_i/i>c_{n-1}^{-a/2}/2$.

\end{lemma}

\begin{pf}

This follows from condition LC2 and Lemma \ref {lands cool}.
\end{pf}

Using that $v_n>c_{n-1}^{-a}$ and $k_i>v_n$ we get

\begin{cor} \label {high time}

With total probability, for $n$ big enough and for
$1 \leq i \leq s_n$,
$$
\frac {\ln(k_i)} {\ln(c_{n-1}^{-1})} >
a/3 \left( 1+\frac {\ln(i)} {\ln(c_{n-1}^{-1})} \right).
$$

\end{cor}

\subsubsection{Upper and lower bounds}

Considering $|R_n(0)|=|f^{v_n}(0)|<c_{n-1}$ and using $v_n<c_{n-1}^{-b}$
we get
$$
\limsup_{n \to \infty} \frac {-\ln |f^n(0)|} {\ln (n)} > a.
$$

Let now $v_n \leq k<v_{n+1}$.  If $|f^k(0)|<k^{-3b^3}$ we have
$f^k(0) \in I_n$ and so $k=k_i$ for some $i$.
It follows from Corollaries \ref {distance} and \ref {high time} that
$$
|f^k(0)|>k^{-3b^3}.
$$

\end{document}